\documentclass{article}
\usepackage{amsfonts}
\usepackage{graphicx}
\usepackage{epstopdf}
\usepackage{amssymb,amsmath,amscd}
\usepackage{mathrsfs}
\usepackage{amsfonts}
\date{}
\newcommand\be{\begin{array}}
\newcommand\en{\end{array}}
\newcommand\di{\displaystyle}

\newcommand\ga{\gamma}
\newcommand\ep{\emptyset}

\newcommand\si{\sigma}

\newcommand\na{\nabla}
\newcommand\vr{\varphi}

\newcommand\pa{\partial}

\newcommand\va{\varepsilon}
\newcommand\la{\lambda}
\newcommand\Om{\Omega}
\newcommand\wi{\widetilde}
\newcommand\de{\delta}
\newcommand\De{\Delta}

\newcommand\gl{\geqslant}
\newcommand\ls{\leqslant}

\newcommand\ri{\rightarrow}
\newcommand\ti{\times}

\newcommand\bk{\backslash}
\newcommand\al{\alpha}
\newcommand\lge{\langle}
\newcommand\rge{\rangle}

\newcommand\no{\eqno}
\newcommand\mb{\mathbb}

\begin{document}
\title{\bf     The Method of Invariant Sets of Descending Flow for Locally Lipschitz  Functionals
\thanks{This paper is supported by NSFC 11871250.}
}
\author{
Xian  Xu$^1$,  Baoxia Qin$^2$} \maketitle

\vspace{2mm} \noindent \begin{center}\ $^1$Department of Mathematics,
Jiangsu Normal University, Xuzhou,
\\Jiangsu,
221116,  P. R. China

$^2$School of Mathematics,
Qilu Normal  University, Jinan,
\\Shandong,
250013,  P. R. China\\

\end{center}
\begin{center}
\begin{minipage}{5in}
{\small {\bf Abstract}\quad  In this paper, we extend the method of invariant sets of descending flow that proposed by Sun Jingxian for smooth functionals to the locally Lipschitz functionals. By  this way, we obtain the existence results for the positive, negative and sign-changing critical points of the  locally Lipschitz functionals, and apply these theoretical results to the study of differential inclusion problems with $p$-Laplacian.  In order to obtain the above results, we develop some new techniques: 1) We establish the method of how to extend the pseudo-gradient field to the whole space on the premise of preserving the useful information of the local pseudo-gradient field; 2) In the case of set-valued mapping, a pseudo-gradient field is established to make both the cone and the negative cone being invariant sets of descending flow. To obtain our main results, a new class of (PS) condition is also proposed.

}
\end{minipage}
\end{center}

\section{Introduction}

\ \ \ \ \    \rm The method of invariant sets of descending flow was proposed by Sun Jingxian in [2]. So far, it has been widely used to study the existence of solutions of elliptic equation boundary value problems.
Now let us recall some theories of the method of invariant sets of descending flow.
Let~$E$ be a real~Banach space, $J:E\rightarrow \mathbb{R}$  a~$C^1$ functional. Let~$K=\{u\in E: J'(u)=0\}, E_0=E\setminus K, W:E_0\rightarrow E$ be a pseudo-gradient vector field of~$J$. Now we consider the following initial value problem in~$E_0$
$$
\left\{
\begin{array}{ll}
\displaystyle \frac{d u(t)}{dt}=-W\left(u(t)\right)~~t\gl 0,\\
u(0)=u_0.
\end{array}
\right.
\no(1.1)$$
By the theory of ordinary differential equations in Banach spaces, (1.1) has a unique solution, denoted by $u(t,u_0)$, with its right maximal existence interval $[0, T(u_0))$.
\vskip 0.1in

The following Definition 2.1,2.2, and Theorem 1.1,1.2 can be fund in [1,2].

{\bf Definition 1.1.}(Invariant sets of descending flow) \textit{Let~$M$ be a nonempty subset of~$E$, we call~$M$ is a invariant set of descending flow generated by~$W$ of~$J$, if~$\{u(t, u_0):0\ls t< T(u_0)\}\subset M$ hold for all~$u_0\in M\setminus K$.}
\vskip 0.1in

{\bf Theorem 1.1.}\
\textit{If~$J\in C^1(E, \mathbb{R}), M\subset E$ is a closed invariant set of descending flow of~(1.1), let~$\alpha= \inf\limits_{u\in M} J(u)>-\infty$, $J$ satisfies the~$(PS)$ condition on~$M$, then~$\alpha$ is a critical value of~$J$, and there exist~$u_0\in M$ such that~$J'(u_0)=0, J(u_0)=\alpha$.}
\vskip 0.1in

Theorem~1.1 states that finding different critical points may  be attributed to finding different invariant sets of descending flow.
\vskip 0.1in

{\bf Definition 1.2.}\
\textit{Let~$M$ and~$D$ be invariant sets of descending flow of~$J$, $D\subset M$, let
$$C_M(D)=\{u_0\in M: \mbox{there exists} ~0\ls t'< T(u_0)~\mbox{such~that} ~u(t', u_0)\in D\}.$$
Then~$C_M(D)$ is a invariant set of descending flow of~$J$ expanded by~$D$. If~$D=C_M(D)$, then~$D$ is called a complete invariant set of descending flow of~$J$ with respect to~$M$.}
\vskip 0.1in

{\bf  Theorem 1.2.}\
\textit{Let~$M$ be a closed and connected invariant set of descending flow of~$J$, $D\subset M$ is an open invariant sets of descending flow of~$J$. If~$C_M(D)\neq M, \inf\limits_{u\in {\partial_M} D} J(u)>-\infty$ and~$J$ satisfies the~$(PS)$ condition on~$M\setminus D$, then
$$c=\inf\limits_{u\in \partial_M C_M(D)} J(u)\gl \inf\limits_{u\in {\partial_M} D} J(u)>-\infty,$$
and $c$ is a critical value of~$J$, there exist critical points~$u^*\in \partial_M C_M(D)$ corresponding to~$c$ of~$J$.}
\vskip 0.1in

 Theorem 1.2 states that  there would exist a new critical point if $\partial_M C_M(D)\cap \pa_M (D)=\ep$.  In order to determine whether a closed convex set is an invariant set of descending flow, Sun [3] also introduced the Schauder invariance condition. For more detailed theories on Sun's method, one can refer to the excellent paper  [1] by Liu and  Sun. This method has been extensively  applied by many authors in the past 20 years or so to study the existence of solutions for various of boundary value problems; see [1,4-8,31] and the references therein. For instance,  Liu and Sun [1] applied the method to semilinear elliptic boundary value problems and to second-order Hamiltonian systems. A typical result in [1] is as follows. Consider the semilinear Dirichlet problem
$$ -\Delta u=f(x,u)\ \mbox{in}\ \Omega, \ u=0\ \mbox{on}\ \pa\Om, \no(1.2)$$ where $\Omega\subset \mb R^N$ is a bounded smooth domain and $f\in C(\bar\Omega\ti \mb R)$ is a superlinear function with subcritical growth. Assume that $f(x,u)+mu$ is increasing for some $m>0$. If (1.2) has a supersolution $\phi$ and a subsolution $\psi$ such that $\phi\ls \psi$ and both $\phi$ and $\psi$ are not solutions of (1.2), then (1.2) has at least four solutions.

Chang [9], in ordered  to study  equations with discontinuities, developed an extension of the classical smooth critical point theory, to non-smooth locally Lipschitz functionals.   The theory of Chang was based on the sub-differential of  locally Lipschitz functionals due to Clarke  [10]. Using this sub-differential, Chang proposed a generalization of the deformation lemma and obtained a  Mountain Pass Theorem for Locally Lipschitz functionals. Subsequently, many other saddle points theorems for locally Lipschitz functionals have been obtained and been applied to  the studying of the existence of solutions for various of  elliptic equations boundary value with discontinuities; see [11-24] and the references therein.

 Now, a natural question  arises is \textit{whether the method of invariant sets of descending flow  proposed by Sun can be applied to locally Lipschitz functionals}. This paper will serve to fulfill this purpose. We will try to use this method to obtain  the existence results for  critical points of locally Lipschitz functionals.  A direct extension of this method to the  locally Lipschitz functional is faced with serious difficulties. As is well known, to use the method of invariant sets of descending flow for $C^1$ functionals we need first construct a pseudogradient vector field over the  Banach space. However, for the Lipschitz functional case, one usually only can construct locally a pseudogradient vector field over a subset of the Banach space rather than the whole Banach space. This is the first difficulty need to overcome. In this paper, we will analyze the altitude of the critical point of the locally Lipschitz functionals in advance, construct locally a  pseudogradient vector field on a neighborhood of the  altitude, and then extend this pseudogradient vector field to the entire Banach space. We can then use the method of invariant sets of descending flow that Sun have proposed to obtain the existence results for the critical points for locally Lipschitz functionals. In the course of the study, we need to determine whether a closed convex set is a  invariant set of the descending flow generated by the pseudogradient vector field. In the cases of the functionals being of $C^1$,  we use the Schauder invariance condition presented in the literature [1-3]. However, in the cases of the  functionals being of locally Lipschitz, no one has yet used the Schauder invariance condition to establish invariant flows on convex closed sets. This is the second difficulty we face. To this end, based on the conclusion in  [20] about the relationship between the critical points on the whole space and the critical points on the closed convex sets under the Schauder invariance condition, and using the method in [9,14] to establish the   descending flow on  closed convex sets, we get the result that the closed convex set is a  descending flow invariant set. Using the above method, in this paper we get the existence results   of at least one sign-changing, at least one positive and at least one negative critical point of the locally Lipschitz functionals. The theoretical results can be applied to the  study of the existence of sign-changing solutions for differential inclusion problems with $p$-Laplacian and so extend the relevant results in literatures.

\section{ Main Results}
\ \ \

In what follows of this section we will let $X$ and $E$ be two  real Banach spaces with the norms $\|\cdot\|$ and $\|\cdot\|_1$, respectively.  Assume that  $X$ is reflexive, $E$ is densely imbedded in $X$.   Let $X^*$ be the topological dual of $X$ and $\lge\cdot, \cdot\rge$ denote the duality pairing between $X^*$ and $X$.  Let $P$ be a closed convex cone of $X$, that is, $P$ is closed convex set in $X$, $\la x \in P$ for all $x\in P$ and $\la\gl 0$, and $P\cap (-P)=\{0\}$. Let $P_1= P\cap E$. We assume that $P_1$ has a nonempty interior in the $E$ topology, and denote its interior in the $E$ topology by int$P_1$. For $x\in X$ and $A\subset X$, let $\mbox{dist}_X(x,A)=\inf\limits_{y\in A}\|x-y\|$. For any $R>0$, let $B(0, R)=\{x\in X: \|x\|<R\}$, $B_1(0, R)=\{x\in E: \|x\|_1<R\}$ and $S_R=\{x\in X: \|x\|=R\}$. Given a subset $A\subset E$, we write $\pa_E A$  for the boundary of $A$ in $E$.

Let us recall some theories concerning the  sub-differential theory of locally Lipschitz functionals due to Clarke [10].   A functional $\vr :X\ri \mb R$ is said to be locally Lipschitz, if for every $x\in X$, there exists a neighbourhood $U$ of $x$ and a constant $k>0$ depending on $U$ such that $|\vr (z)-\vr(y)|\ls k \|z-y\|$ for all $z,y\in U$.  For such a  functional we define generalized directional derivative $\vr^0(x; h)$ at $x\in X$ in the direction $h\in X$ by
$$\vr ^0(x;h)=\lim\limits_{x'\ri x}\sup\limits_{\la\downarrow 0^+}\frac{\vr(x'+\la h)-\vr (x')}{\la}.$$
The function $h\mapsto \vr^0(x; h)$ is sublinear, continuous. So by the Hahn-Banach theorem we know that $\vr^0(x; \cdot)$ is the support function of a nonempty, convex and $w^*$-compact set
$$\pa\vr (x)=\{x^*\in X^*:\lge x^*,h\rge\ls \vr^0(x;h)\ \mbox{for all}\ h\in X\}.$$
The set $\pa\vr(x)$ is called the generalized or Clarke sub-differential of $\vr$ at $x$.  A point $x\in X$ is a critical point of  $\vr$ if $0\in \pa\vr(x)$. Let $\mb K=\{x\in X: 0\in\pa \vr(x)\}$.
\vskip 0.1in

{\bf Proposition 2.1 ([12,13]).}\ \textit{1)\ If $\vr,\psi:X\ri\mb R$ are locally Lipschitz functionals, then $\pa(\vr+\psi)(x)\subset \pa\vr (x)+\pa\psi(x)$, while for any $\la\in \mb R$ we have $\pa(\la\vr)(x)=\la\pa\vr(x)$; 2)\ If $\vr:X\ri \mb R$ is also convex, then this sudifferential coincides with the sub-differential in the sense of convex analysis. If $\vr$ is strictly differentiable, then $\pa\vr(x)=\{\vr'(x)\}$; 3)\ If $\vr:X\ri\mb R$ is locally Lipschitz functional, $\pa \vr(u)$ is a weakly$^*$-compact subset of $X^*$ which is bounded by the Lipschitz constant $k>0$ of $\vr$ near $u$.  }

\vskip 0.1in

S.T. Kyritsi and N. S. Papageorgiou  [14] developed a  critical  point theory for nonsmooth locally Lipschitz functionals defined  on closed, convex set  extending this way the work of Struwe.  Let $C\subset X$ be a nonempty, non-singleton, closed and convex set.  For $x\in C$ we define
 $$m_C(x)=\inf\limits_{x^*} \sup\limits_{y}\big\{\lge x^*,x-y\rge:y\in C,\|x-y\|<1,x^*\in\pa \vr(x)\big\}.$$
 Evidently, $m_C(x)\gl 0$ for all $x\in C$. This quantity can be viewed as a  measure of the generalized slope of $\vr$ at $x\in C$. If $\vr$ admits an extension $\hat\vr\in C^1(X)$, then  $\pa\vr(x)=\{\vr'(x)\}$ and so we have
 $$m_C(x)=\sup\big\{\lge \vr'(x),x-y\rge:y\in C,\|x-y\|<1\big\},$$
 which is the quantity used by Struwe [29,p.147]. Also if $C=X$, then we have
 $$m_C(x)=m(x)=\inf\{\|x^*\|_*: x^*\in \pa\vr(x)\},$$
 which is the quantity used by Chang [9].
\vskip 0.1in

Now let us introduce the outwardly directed condition and the Schauder invariance condition for set value mappings in a manner as [20]. Let $X$ be  reflexive. As usual, we will identify $X^{**}$ with $X$ while $\textsf{ F}:X^*\mapsto 2^X$  will denote the duality map, given by
 $$\textsf{F}(x^*):=\{x\in X:\lge x^*,x\rge=\|x^*\|^2_{*}=\|x\|^2\},\ \  \forall x^*\in X^*.$$
The set $\textsf{F}(x^*)$ turns out to be nonempty, convex, and closed; see,e.g. [13, pp. 311-319]. As in [20], we define
$$\nabla \vr(x):=\textsf{F}(\pa\vr(x)),\ \  x\in X. \no(2.1)$$
Clearly, $\nabla \vr(x)$ depends on the choice of the duality pairing between $X$ and $X^*$ whenever it is compatible with the topology of $X$. If $X$ is a Hilbert space, the duality paring becomes the scalar product and (2.1) gives the usual gradient. Write $I$ for the identity operator on $X$.
\vskip 0.1in

Suppose $C$ is a convex and closed set of $X$. Let $\de_C:X\ri \mb R\cup \{+\infty\}$ be the indicator function of $C$, namely
 $$\de_C(x):=\left\{\be{ll} 0, &\mbox{if}\ x\in C,\\
 +\infty,&\mbox{otherwise}.
 \en\right.
 $$
 Then we have
 $$\pa\de_C(x)=\big\{x^*\in X^*:\lge x^*,z-x\rge\ls 0, \forall z\in C\big\}.$$
The  set $\pa\de_C(x)$ is usually called normal cone to $C$ at $x$.
\vskip 0.1in

 {\bf Definition 2.1 ([20]).}\  \textit{ Suppose $X$ is reflexive, $C$ is a convex and closed set of $X$, $\vr:X\mapsto \mb R$ is  locally Lipschitz continuous. If for  $x\in \pa C$,
$$(-\pa\de_C(x))\cap \pa \vr(x)\subset\{0\},$$
then  we say that $\pa \vr$ turns out to be outwardly directed on $ \pa C$. This clearly rewrites as
$$\forall z^*\in \pa\vr(x)\backslash\{0\}\ \mbox{there exists}\ z\in C\ \mbox{fulfilling}\ \lge z^*, z-x\rge<0.$$}
 \vskip 0.1in

 {\bf Definition 2.2 ([20], Schauder invariance condition).}   \textit{Suppose $X$ is reflexive, $C$ is a convex and closed set of $X$, $\vr:X\mapsto \mb R$ is  locally Lipschitz continuous. Then we call $\vr$ satisfies the Schauder invariance condition on $\pa C$ if
 $\big(I-\na\vr\big)(\pa C)\subset C$.}
  \vskip 0.1in

  {\bf Remark 2.1.}\ The well known Schauder invariance condition for a $C^1$-functional $\vr$ on a Hilbert space $X$ reads as $(I-\vr')(C)\subset C$; see [1-3,25]. It has been extended to  Banach spaces in [26]. The notation of Schauder invariance condition was firstly put forward by Sun Jingxian in [3].   It follows from [20, Theorem 4.4 and 4.5] we have the following Lemma 2.1.
\vskip 0.1in

 {\bf Lemma 2.1 ([20]).}\   \textit{Suppose $X$ is reflexive, $C$ is a convex and closed set of $X$, $\vr:X\mapsto \mb R$ is  locally Lipschitz continuous,  $\vr$ is outwardly directed at $\pa C$ or $\vr$ satisfies the Schauder invariance condition on $\pa C$. Then  $m_C(x)=0$ if and only if $m(x)=0$.}
 \vskip 0.1in

{\bf Lemma 2.2 ([36],  Von Neumann).}\  \textit{ Let $X, Y$ be two Hausdorff topological linear spaces, $C\subset X$, $D\subset Y$ be two convex and compact sets. Let $\psi: X\ti Y\mapsto \mb R$ satisfy:
1)\ $x\mapsto \psi(x, y)$ is upper semi-continuous and concave; 2)\ $y\mapsto \psi(x, y)$ is lower semi-continuous and convex.
Then $\psi$ has at least one saddle point $(\bar x, \bar y)\in C\ti D$.}
\vskip 0.1in

{\bf Definition 2.3.}\  \textit{Let $D$ be a nonempty closed subset of $X$. We say that $\vr$ satisfies  the non-smooth $CPS$-condition on $D$, denoted by $(CPS)_{D}$, if every sequence $\{x_n\}\subset D$ such that $\{\vr(x_n)\}$ is bounded and $(1+\|x_n\|) m_D(x_n)\ri 0$ as $n\ri \infty$, has a convergent subsequence. If $D=X$, denoted simply by $(CPS)$. }
\vskip 0.1in

 Let us introduce the following conditions:
  \vskip 0.1in

 (H$_1)$\ $\vr: X\mapsto \mb R$ satisfies the conditions  $(CPS)$, $(CPS)_P$ and $(CPS)_{-P}$;

 (H$_2)$\  $\vr(0)=0$ and $0$ is a strictly local minimum point;

 (H$_3)$\  There exist sub-space $E_1$ of $E$ and $R_0>0$ such that $\mbox{dim}\ E_1=2$, $E_1\cap \mbox{int}(P_1)\neq \ep$ and
 $$\sup\limits_{u\in S_{R_0}\cap E_1} \vr(u)< 0;\no(2.2)$$

(H$_4)$\ either $\vr$ is outwardly directed at $\pa (\pm P)$, or $(I-\na \vr)(\pm P)\subset \pm P$;

  (H$_5)$\  $\big(\mb K\bk \{0\}\big)\cap \big( \pa_E P_1\cup \pa_E (-P_1)\big)=\ep$;

(H$_6)$\ for each $a<b$, $\vr^{-1}([a,b])\cap K$ is compact in $E$.
 \vskip 0.1in

 {\bf Remark 2.2.}\ It follows from [20, Theorem 4.5] that $\vr$ is outwardly directed at $\pa (\pm P)$ if $(I-\na \vr)(\pm P)\subset \pm P$. Here, in condition (H$_4)$,  we list the above two conditions at the same time for the purpose of application convenience.
\vskip 0.1in

  We have the following main result.

{\bf Theorem 2.1.}\ \textit{Suppose that  (H$_1)\sim $(H$_6)$ hold. Then $\vr$ has at least one positive, one negative and one sign-changing critical point.}
  \vskip 0.1in

It follows from $E_1\cap\mbox{int}\ P_1\neq\ep$ that $E_1\cap \mbox{int}(-P_1)\neq\ep$. Let $d_1=\sup\limits_{x\in B(0, {R_0})\cap E_1} \vr(x)+1$ and $D_0=\vr^{-1}\big([-1, d_1]\big)$.  Let $K=\mb K\cap D_0$.   Because $0$ is an isolate critical point, we may take $\bar R_1>0$ small enough such that $B(0,\bar R_1)\cap (K\bk\{0\})=\ep$.  Since $E\hookrightarrow X$, we take $R_1>0$ small such that $B_1(0, R_1)\subset B(0, \bar R_1)$. By using the  condition (H$_1)$,(H$_5)$ and (H$_6)$ we  may take $\de>0$ small enough such that $$ K_{3\de}\cap \big( \pa_E P_1\cup \pa_E (-P_1)\cup B_1(0, R_1)\big)=\ep,\no(2.3)$$
where $K_{3\de}=\{x\in D_0: \mbox{dist}_E (x, K)<3\de\}$.

To show Theorem 2.1 we need to give some lemmas.
 \vskip 0.1in

 {\bf Lemma 2.3.}\ \textit{There exists a locally Lipschitz mapping $v: D_0\bk K_\de\mapsto X$ such that $\|v(x)\|\ls 2(1+\|x\|)$ for any $x\in D_0\bk K_\de$, $\lge x^*, v(x)\rge\gl\frac{\ga}{16}$ for some $\ga>0$ and all $x^*\in\pa \vr(x)$. Moreover, $v:D_0\bk K_\de\mapsto E$ is locally Lipschitz and
$$x- \frac{1}{1+\|x\|} v(x)\in P_1\ \mbox{for any }\ x\in P\cap (D_0\bk K_\de), \no(2.4) $$
 $$x- \frac{1}{1+\|x\|} v(x)\in -P_1\ \mbox{for any }\ x\in -P\cap (D_0\bk K_\de).\no(2.5)$$}

 {\bf Proof.}\ Let $S=D_0\bk(P\cup (-P))$. First we claim that
$$(1+\|x\|)m(x)\gl \ga,\ \ \  \forall x\in S\bk K_\de,\no(2.6)$$
$$(1+\|x\|)m_P(x)\gl \ga,\ \ \  \forall x\in  (D_0\cap P)\bk K_\de, \no(2.7)$$
$$(1+\|x\|)m_{-P}(x)\gl \ga,\ \ \  \forall x\in (D_0\cap (-P))\bk K_\de. \no(2.8)$$

We only show that (2.7) holds. In a similar way we can show that (2.6) and (2.8) hold. Arguing by make contradiction that (2.7)  does'n  hold. Then  there exists an sequence $\{x_n\}\subset (D_0\cap P)\bk K_\de$ such that $(1+\|x_n\|)m_P(x_n)\ri 0$ as $n\ri \infty$. Obviously, $\{\vr(x_n)\}$ is bounded. Since $\vr$ satisfies the $(CPS)_P$ condition, up to a subsequence if necessary, we may assume that $x_n\ri x_0$ as $n\ri\infty$ for some $x_0\in (D_0\cap P)\bk K_\de$. Since $m_P(\cdot):P\mapsto \mb R$ is lower semi-continuous, we have $m_P(x_0)=0$. It follows from  (H$_4)$ and Lemma 2.1 that $m(x_0)=0$, which is a contradiction. Thus, (2.7) holds.

For each $x_0\in S\bk K_\de$, take $w_0^*\in \pa \vr(x_0)$ such that $m(x_0)=\|w_0^*\|_{*}>0$. Then we have $B(0, \|w_0^*\|_{*})\cap \pa \vr(x_0)=\ep$, where $$B(0, \|w_0^*\|_{*})=\{z^*\in X^*: \|z^*\|_{*}<\|w_0^*\|_{*}\}.$$ So,  by using the separation theorem in the weak$^*$-topology, we can find $u_1(x_0)\in X$ with $\|u_1(x_0)\|=1$ such that for all $z^*\in B(0, \|w_0^*\|_{*})$ and $y^*\in\pa \vr(x_0)$,
$$\lge z^*, u_1(x_0)\rge \ls \lge w_0^*, u_1(x_0)\rge\ls \lge y^*, u_1(x_0)\rge.$$
Recall that $$\|w_0^*\|_{*}=\sup\big\{\lge z^*, u_1(x_0)\rge: z^*\in B(0, \|w_0^*\|_{*})\big\},$$
then we have by (2.6),
$$\lge y^*, u_1(x_0)\rge \gl \|w_0^*\|_{*}>\frac{\ga}{2(1+\|x_0\|)}.$$
Since the map $x\mapsto\pa \vr(x)$ is usc. from $X$ into $X^*_w$, we may take an open neighborhood $B_1(x_0, r_1(x_0))$ of $x_0$,  such that
$$\lge  y^*, u_1(x_0)\rge>\frac{\ga}{4(1+\|y\|)},\ \ \forall y^*\in\pa \vr(y),\ y\in U_1(x_0), \no(2.9)$$
where $U_1(x_0)=B_1(x_0, r_1(x_0))\cap (S\bk  K_\de)$.
Since $x_0\in S\bk  K_\de$ and $ S\bk K_\de$ is an open subset of $D_0\bk K_\de$, we may take $r_1(x_0)>0$ small enough such that $U_1(x_0)\subset S\bk K_\de $.

 Pick $x_0\in (P\cap D_0)\bk K_\de$. Let $C=(\{x_0\}-P)\cap \bar B(0,1)$ and $D=\pa\vr(x_0)$, where $\bar B(0, 1)=\{x\in X: \|x\|\ls 1\}$. Let $X_w$ and $X^*_w$ be the spaces $X$ and $X^*$ furnished their weak topology respectively. Since $X$ is reflexive,  $D$ is  compact in $X^*_w$ and $C$ is  compact  in $X_w$. Obviously, both $C$ and $D$ are convex. Let $\psi: C\ti D\mapsto \mb R$ be defined by $\psi(x, y^*)=\lge y^*, x\rge$ for any $(x, y^*)\in C\ti D$. It follows from Lemma 2.2 that $\psi$ has at least one saddle point $(u_2(x_0), x_0^*)\in C\ti D$, that is
$$\min\limits_{x^*\in D}\max\limits_{x\in C}\lge x^*, x\rge=\lge x_0^*, u_2(x_0)\rge =\max\limits_{x\in C}\min\limits_{x^*\in D}\lge x^*, x\rge, \ \forall x^*\in \pa\vr(x_0), x\in (\{x_0\}-P)\cap \bar B(0,1).$$ Thus, we have
$$\lge x_0^*, x\rge\ls\lge x_0^*, u_2(x_0)\rge\ls \lge x^*, u_2(x_0)\rge, \ \forall x^*\in \pa\vr(x_0), x\in (\{x_0\}-P)\cap \bar B(0,1).$$   It follows from (2.7) that for any $x^*\in\pa\vr(x_0)$,
$$\lge x^*, u_2(x_0)\rge\gl \lge x_0^*, u_2(x_0)\rge=m_{P}(x_0)>\frac{\ga}{2(1+\|x_0\|)}.$$
By using the fact that $x\mapsto\pa \vr(x)$ is upper semi-continuous, we know that there exists an open neighborhood $B_2(x_0, r_2(x_0))$ of $x_0$  such that for any $y\in U_2(x_0)$, $y^*\in\pa \vr(y)$, $$\lge y^*,u_2(x_0)\rge>\frac{\ga}{4(1+\|y\|)}.\no(2.10)$$
where $U_2(x_0)=B_2(x_0, r_2(x_0))\cap (D_0\bk K_\de)$. Since $x_0\in (P\cap D_0)\bk K_\de\subset D_0\bk \big(K_\de\cup (-P)\big)$ and $D_0\bk \big(K_\de\cup (-P)\big)$ is an open set of $D_0\bk K_\de$, we may take $r_2(x_0)>0$ small  such that  $U_2(x_0)\cap (-P)=\ep$.

Similarly, by (2.8) we can show that for each $x_0\in (D_0\cap (-P))\bk K_\de$, there exist $r_3(x_0)>0$, $u_3(x_0)\in X$ with $\|u_3(x_0)\|\ls 1$, such that $$x_0-u_3(x_0)\in -P,\ \  U_3(x_0)\cap P=\ep$$ and
for any $y\in U_3(x_0)$, $y^*\in\pa \vr(y)$, $$\lge y^*,u_3(x_0)\rge>\frac{\ga}{4(1+\|y\|)},\no(2.11)$$
where $U_3(x_0):=B_3(x_0, r_3(x_0))\cap (D_0\bk K_\de)$.

By 3) in Proposition 2.1, we may assume that $\|x^*\|_{*}\ls L_{\al, i}$ for some $L_{\al, i}>0$ and any  $x\in U_i(x_\al) $ with $i=1,2,3$ and $x^*\in\pa \vr(x)$.  Also, we assume that for $i=1,2,3$, $B_i(x_\al, r_i(x_\al))$ has a small radium $r_i(x_\al)>0$ such that $(1+\|x\|)(1+\|x_\al\|)^{-1}\ls 2$ for each $x\in U_i(x_\al)$, and $$0< r_i(x_\al)\ls \min\Big\{\frac{1}{2},\frac{\ga}{16(1+\|x_\al\|) L_{\al, i}}\Big\}.$$

Let $$\mathscr A_1=\big\{U_1(x_0): x_0\in S\bk K_\de\big\},$$
$$\mathscr A_2= \big\{U_2(x_0): x_0\in (D_0\cap P)\bk K_\de\big\}, $$
$$\mathscr A_3= \big\{U_3(x_0): x_0\in (D_0\cap (-P))\bk K_\de\big\},$$
 and $\mathscr A=\mathscr A_1\cup \mathscr A_2\cup \mathscr A_3$. Then $\mathscr A$ is an open cover of $D_0\bk K_\de$. By paracompactness we can find a locally finite refinement $\mathscr B=\{V_\al:\al\in\Lambda\}$ and a locally Lipschitz partition of unit $\{\ga_\al:\al\in\Lambda\}$ sub-ordinate to it. For each $\al\in\Lambda$ we can find $x_\al\in D_0\bk K_\de$ such that $V_\al\subset U_{i(\al)}(x_\al)$ for some  $i(\al)\in \{1,2,3\}$,   and $U_{i(\al)}(x_\al)\in \mathscr A$.  To this $x_\al\in D_0\bk K_\de$ corresponds  the element $w_\al^{i(\al)}$ such that $\|w_\al^{i(\al)}\|\ls 1$, and $w_\al^{i(\al)}=u_{i(\al)}(x_\al)$ if $V_\al\subset U_{i(\al)}(x_\al)$ for some $i(\al)\in \{1,2,3\}$.
Since $E$ is densely imbedded in $X$, we may take $\bar x_\al, \bar w_\al^{i(\al)}\in E$, such that $\|\bar w_\al^{i(\al)}\|\ls 1$, $\bar x_\al-\bar w_\al^{i(\al)}\in P_1$ when $x_\al\in P$,  $\bar x_\al-\bar w_\al^{i(\al)}\in -P_1$ when $x_\al\in -P$ and
$$\max\big\{\|w_\al^{i(\al)}-\bar w_\al^{i(\al)}\|, \|x_\al-\bar x_\al\|\big\}<\min\Big\{\frac{1}{2},\frac{\ga}{32(1+\|x_\al\|) L_{\al, i(\al)}}\Big\}.\no(2.12)$$

Now, let $v:D_0\bk K_\de\mapsto X$ be defined by
$$v(x)=(1+\|x\|)\sum\limits_{\al\in \Lambda} \ga_\al(x) (\bar w_\al^{i(\al)}-\bar x_\al+x).\no(2.13)$$
It is easy to see that $v:D_0\bk K_\de\mapsto X$ is locally Lipschitz. Since $E\hookrightarrow X$, by [32, Lemma 2.3] we see that $\ga_\al: E\to\mb R$ is also locally Lipschitz. So, we can prove that  $v:D_0\bk K_\de\mapsto E$ is locally Lipschitz.

By (2.12) and (2.13),  we have $$\|v(x)\|\ls(1+\|x\|)\sum\limits_{\al\in \Lambda} \ga_\al(x) (\|\bar w_\al^{i(\al)}\|+\|\bar x_\al-x_\al\|+\|x_\al-x\|)\ls 2(1+\|x\|),$$
and
$$\be{ll}(1+\|x\|)&\big| \sum\limits_{\al\in \Lambda} \ga_\al(x) \lge x^*, x-\bar x_\al\rge\big|\\
&\ls \sum\limits_{\al\in\Lambda}\frac{1+\|x\|}{1+\|x_\al\|}(1+\|x_\al\|)\ga_\al(x)\|x^*\|_{*}\big(\|x- x_\al\|+\|x_\al-\bar x_\al\|\big)\\
&\ls 2 \sum\limits_{\al\in\Lambda} L_{\al, i(\al)}\big(\|x- x_\al\|+\|x_\al-\bar x_\al\|\big)(1+\|x_\al\|)\\
&<\frac{\ga}{16}\en
$$
for each $x\in D_0\bk K_\de$.  On the other hand,$$\be{ll}\Big|\sum\limits_{\al\in\Lambda}\ga_\al(x)(1+\|x\|)\lge x^*, \bar w_\al^{i(\al)}-w_\al^{i(\al)}\rge\Big|&\ls \sum\limits_{\al\in\Lambda}\frac{1+\|x\|}{1+\|x_\al\|}(1+\|x_\al\|)\ga_\al(x)\|x^*\|_{*}\|\bar w_\al^{i(\al)}-w_\al^{i(\al)}\|\\
&\ls 2\sum\limits_{\al\in\Lambda}(1+\|x_\al\|)L_{\al, i(\al)}\|\bar w_\al^{i(\al)}-w_\al^{i(\al)}\|<\frac{\ga}{16},\en $$
and so, by (2.9)$\sim$(2.11) we have $$\be{ll}\sum\limits_{\al\in\Lambda}\ga_\al(x)(1+\|x\|)\lge x^*, \bar w_\al^{i(\al)}\rge&=\sum\limits_{\al\in\Lambda}\ga_\al(x)(1+\|x\|)\lge x^*, \bar w_\al^{i(\al)}-w_\al^{i(\al)}+w_\al^{i(\al)}\rge\\
&= \sum\limits_{\al\in\Lambda}\ga_\al(x)(1+\|x\|)\lge x^*, w_\al^{i(\al)}\rge\\
&+\sum\limits_{\al\in\Lambda}\ga_\al(x)(1+\|x\|)\lge x^*, \bar w_\al^{i(\al)}-w_\al^{i(\al)}\rge\\
&\gl \frac{\ga}{4}-\frac{\ga}{16}\gl\frac{\ga}{8}.\en$$
Then,  we have
$$\be{ll} \lge x^*, v(x)\rge&= \sum\limits_{\al\in\Lambda}\ga_\al(x)(1+\|x\|)\lge x^*, \bar w_\al^{i(\al)}\rge\\
&+\sum\limits_{\al\in\Lambda}\ga_\al(x)(1+\|x\|)\lge x^*, x-\bar x_\al\rge\gl \frac{\ga}{16}.\en$$

For any $x\in P\cap (D_0\bk K_\de)$, we have
$$x-\frac{1}{1+\|x\|} v(x)=\sum\limits_{\al\in\Lambda}\ga_\al(x)(\bar x_\al- \bar w_\al^{i(\al)})\in P_1.$$
This implies that (2.4) holds. Similarly, (2.5) also holds. The proof of  is complete.
\vskip 0.1in

Let $l_1, l_2:D_0\mapsto \mb R$ be defined by  for any $ x\in D_0$,
$$l_1(x)=\frac{\mbox{dist}_X\ (x,  K_{\de})}{\mbox{dist}_X\ (x, K_\de)+ \mbox{dist}_X\ (x, D_0\bk K_{2\de})},$$
$$l_2(x)=\frac{\mbox{dist}_X\ \big(x,  \big(\vr^{-1}([d_1-\frac{1}{4}, d_1])\cup \vr^{-1}([-1, -\frac{1}{2}])\big)\big)}{\mbox{dist}_X\ \big(x, \vr^{-1}([-\frac{1}{4}, d_1-\frac{1}{2}])\big)+\mbox{dist}_X\ \big(x,  \big(\vr^{-1}([d_1-\frac{1}{4}, d_1])\cup \vr^{-1}([-1, -\frac{1}{4}])\big)\big)},$$
$$v_1(x)=\left\{\be{ll}v(x), & \mbox{if}\ x\in D_0\bk K_\de;\\
0, &{K_\de},\en\right.
$$
and
$$V(x)=\left\{\be{ll}l_1(x)l_2(x)v_1(x), & \mbox{if}\ x\in D_0;\\
0, &\mbox{othewise.}\en\right.
$$ Then $l_1, l_2: X\mapsto \mb R$  are locally Lipschitz. Since $E\hookrightarrow X$, $l_1, l_2: E\mapsto \mb R$  are also locally Lipschitz. Obviously, $V: X\mapsto X$ is locally Lipschitz.

Consider the following initial value problem
$$\left\{\be{l}\frac{du}{dt}=-V(u),\\ u(0)=v_0\in X.\en\right.\no(2.14)
$$
By the theories for initial value problems of ordinary equations in Banach space, we see that (2.14) has a unique solution $\si(t, v_0)$, with its right maximal existence interval  $[0, T(v_0))$ in the $X$ topology, and its right maximal existence interval  $[0, T_1(v_0))$ in the  $E$ topology. Since $E\hookrightarrow X$, we have $T_1(v_0)\ls T(v_0)$.
\vskip 0.1in

As the Definition 1.1 and 1.2, we give the following Definition 2.3 and 2.4.

{\bf Definition 2.4.}\ \textit{A nonempty subset $D$ of $E$ is called an  invariant set of descending flow of (2.14) if $o(v_0)\subset D$ for all $v_0\in D$, where $o(v_0)=\{\si(t,v_0)\subset E:t\in [0, T_1(v_0))\}$.}
\vskip 0.1in
{\bf Definition 2.5.}\ \textit{Let $M\subset E$ be a connected invariant set of descending flow of (2.14),  $D$ be an open subset of $M$ and be an invariant set of descending flow of (2.14). Denote
$$C_{M}(D)=\{v_0: v_0\in D\ \mbox{or}\ v_0\in M\bk D\ \mbox{and there exists}\  t'\in (0,T_1(v_0))\ \mbox{such that}\ \si(t', v_0) \in D\}.$$
If $D=C_{M}(D)$, then $D$ is called  a complete invariant set of descending flow of (2.14) in $M$.}
\vskip 0.1in

{\bf Lemma 2.4.}\  \textit{For each $v_0\in E$, the solution $\si(t, v_0)$ of (2.14) has the following properties:\\
1)\ $T(v_0)=+\infty$;\\
2)\  $T_1(v_0)=+\infty$ if $o(v_0)\subset D_0\bk K_\de$ and $v_0\in (D_0\cap E)\bk K_\de$;\\
3)\ $\vr(\si(t, v_0))$ is non-increasing in $t\in [0,+\infty)$;\\
4)\ $P_1$, $-P_1$, $\mbox{int}\ P_1$ and $\mbox{int}\  (-P_1)$ are all invariant sets of descending flow of (2.14);\\
5)\ For each $v_0\in \vr^{-1}\big((-\infty, d_1-\frac{1}{2}]\big)\cap E$ with $\inf\limits_{u\in o(v_0)} \vr(u)\gl-\frac{1}{4}$,  there exists $\tau(v_0)\gl 0$ such that $\si(\tau(v_0), v_0)\in K_{2\de}$.}

{\bf Proof.}\ 1)\ First we prove that $T(v_0)=+\infty$.  Arguing by make contradiction that $T(v_0)<+\infty$. By (2.14) we have
$$\|\si(t, v_0)-v_0\|\ls\int^t_0\|V(\si(s, v_0))\|ds\ls 2\int^t_0\big(1+\|\si(s, v_0)\|\big)ds.$$
So, we have
$$\be{ll}\frac{1}{2}\|\si(t, v_0)-v_0\|&\ls\di\int^t_0(1+\|\si(s, v_0)\|)ds\\
&\ls\di\int^t_0\|\si(s, v_0)-v_0\|ds+(1+\|v_0\|) t,\en$$
By the well known Gronwall's inequality, we have
$$\be{ll}\frac{1}{2}\|\si(t, v_0)-v_0\|&\ls\di\int^t_0(1+\|v_0\|)e^{t-s}ds+(1+\|v_0\|) t\\
&\ls (1+\|v_0\|)(e^{t}-1)+(1+\|v_0\|) t\\
&\ls (1+\|v_0\|)(t+e^{t}-1)\\
&\ls (1+\|v_0\|)(T(v_0)+e^{T(v_0)}).\en$$
So, we have
$$\|\si(t, v_0)\|\ls 2(1+\|v_0\|)(T(v_0)+e^{T(v_0)})+2\|v_0\|=: M_0(v_0).$$
Take $\{t_n\}\subset [0, T(v_0))$ such that $t_n\ri T^-(v_0)$ and  for $n=1,2,\cdots$, $$|t_n-t_{n-1}|<\frac{1}{2\cdot 2^n\big(1+M_0(v_0)\big)}.$$
Then we have
$$\be{ll} \|\si(t_n, v_0)-\si(t_{n-1}, v_0)\|&\ls\di\int_{t_{n-1}}^{t_n} \|V(\si(s, v_0)\|ds\\
&\ls 2 \di\int_{t_{n-1}}^{t_n} \big(1+\|\si(s, v_0)\|\big)ds\\
&\ls 2\big(1+M_0(v_0)\big)(t_n-t_{n-1})<\frac{1}{2^n}.\en$$
This implies that $\{\si(t_n, v_0)\}$ is a Cauchy sequence. Thus, there exists $\bar u\in X$ such that $\si(t_n, v_0)\ri \bar u$ as $t_n\ri T^-(v_0)$. Then, we can show that $\si(t, v_0)\ri \bar u$ as $t\ri T^-(v_0)$.

Now we consider the initial value problem
$$\left\{\be{l}\frac{du}{dt}=-V(u),\\
u(0)=\bar u,\en\right.\no(2.15)$$
Then (2.15) has a unique solution on $[0, \bar\de)$ for some $\bar\de>0$, and so (2.14) has a unique solution on $[0, T(v_0)+\bar\de)$, which is a contradiction. Thus, we have $T(v_0)=+\infty$.

2)\ Let the operator $A: E\ri X$ be defined by
$$Ax=x-\frac{1}{1+\|x\|}V(x), \ \ \ \forall x\in E.$$
Then we have for any $x\in E$,
$$Ax=\left\{\be{ll} x,\ &\mbox{if}\ x\in\big( E\bk \vr^{-1}(-\frac{1}{2}, d_1-\frac{1}{4})\big)\cup K_\de;\\
x-\frac{1}{1+\|x\|}v(x),\ &\mbox{if}\ x\in \vr^{-1}([-\frac{1}{4}, d_1-\frac{1}{2}])\bk K_{2\de};\\
l(x)\Big(x-\frac{1}{1+\|x\|}v(x)\Big)+\big(1-l(x)\big)x,\ &\mbox{otherwise}.\en
\right\}\in E,$$
where $l(x):=l_1(x)l_2(x)\in (0,1)$, and
$$Ax=\sum\limits_{\al\in \Lambda} \ga_\al(x) (\bar w_\al^{i(\al)}-\bar x_\al), \ \forall x\in D_0\bk K_{2\de}.\no(2.16)$$

 Let
$$\mu(t)=\int^t_0\big(1+\|\si(s, v_0)\|\big)ds \ \mbox{for}\ t\in [0,+\infty).$$
Obviously, $\mu:[0,T_1(v_0))\ri [0,+\infty)$ is increasing. And so, $\mu^{-1}$, the inverse function of $\mu$, exists. Then we have
$$\left\{\be{l}\di\frac{d}{dt}\big(e^{\mu(t)} \si(t, v_0)\big)=e^{\mu(t)}\big(1+\|\si(t, v_0)\|\big) A\si(t, v_0),\\
\si(0, v_0)=v_0.\en\right.
$$
By direct computation, we have
$$\si(t,v_0)=e^{-\mu(t)} v_0+e^{-\mu(t)}\int^t_0e^{\mu(s)}\big(1+\|\si(s, v_0)\|\big)A\si(s, v_0)ds.\no(2.17)$$
where the integral is in the sense of $X$ topology.

By (2.16), one can easily show that for each $T>0$, $\{A(\si(t, v_0)): t\in [0, T]\}$ is contained in a finite-dimensional subspace of $X$.
Now we show that $T_1(v_0)=+\infty$  when  $o(v_0)\subset D_0\bk K_{2\de}$. Arguing by make contradiction that $T_1(v_0)<+\infty$. Take $T> T_1(v_0)$. Note $\{\si(t, v_0):t\in [0, T]\}$ is bounded in $X$. Then, there exists $ M_1(T)>0$ such that for all $s\in [0, T]$, $$\|e^{\mu(s)}\big(1+\|\si(s, v_0)\|\big)A\si(s, v_0)\|_1\ls M_1(T).$$ Let  $\{t_n\}\subset [0, T_1(v_0))$ such that $t_n\ri T_1^-(v_0)$ as $n\ri\infty$.
Note (2.17) also holds in which  the integral is in the sense of $E$ topology for any $t\in [0,T_1(v_0))$. For each $n=2,3, \cdots$, let $\tau_n=\max\{t_n, t_{n-1}\}$ and $\tau_{n-1}=\min\{t_n, t_{n-1}\}$. Then we have
$$\be{ll} \|\si(t_n, v_0)-\si(t_{n-1}, v_0)\|_1&\ls |e^{-\mu(t_n)}-e^{-\mu(t_{n-1})}|\Big( \|v_0\|_1\\
&+\di\int^{t_n}_0e^{\mu(s)}\big(1+\|\si(s, v_0)\|\big)\|A\si(s, v_0)\|_1ds\Big)\\
&+e^{-\mu(t_{n-1})}\di\int^{\tau_n}_{\tau_{n-1}}e^{\mu(s)}\big(1+\|\si(s, v_0)\|\big)\|A\si(s, v_0)\|_1ds\\
&\ls |e^{-\mu(t_n)}-e^{-\mu(t_{n-1})}|\Big( \|v_0\|_1+TM_1(T)\Big)+M_1(T)|t_n-t_{n-1}|.
\en$$
So, $\{\si(t_n, v_0)\}$ is a Cauchy sequence in $E$. Assume that   $\si(t_n, v_0)\ri \bar u$  in $E$ as $t_n\ri T_1^-(v_0)$.
Since $V: E\mapsto E$ is also locally Lipschitz,  one can easily get a contradiction as the proof of $T(v_0)=+\infty$. Thus, we have $T_1(v_0)=+\infty$.

3)\ Let $h(t, v_0)=\vr(\si(t, v_0))$ for all $t\in [0,T_1(v_0))$. It is easy to see that $h(t, v_0)$ is locally Lipschitz in $t\in [0,+\infty)$, hence differentiable almost everywhere. According to Leburng's Mean Theorem  we have
$$\be{ll}\frac{\pa }{\pa s}h(s, v_0)&\ls \max\big\{\lge w^*, \frac{\pa}{\pa s} \si(s, v_0)\rge: w^*\in \pa \vr(\si(s, v_0))\big\}\ \mbox{a.e.}\\
  &=-\min\big\{ \lge w^*, V( \si(s, v_0))\rge: w^*\in \pa \vr(\si(s, v_0))\big\}\ \mbox{ a.e.}\\
 &\ls\left\{\be{ll}-\frac{\ga}{16}, &\mbox{if}\ \si(s, v_0)\in \vr^{-1}([-\frac{1}{4}, d_1-\frac{1}{2}])\bk K_{2\de};\\
 0, &\mbox{otherwise}\en\right.\en\no(2.18)$$
for almost every $t\in [0, T_1(v_0))$. Consequently, $\vr(\si(t,v_0))$ is non-increasing in $t\in [0, T_1(v_0))$.

4)\ For $u\in P_1$, it follows from Lemma 2.3 that for $\la>0$ small enough,
$$ u+\la (-V(u))
=\la l_1(u)l_2(u)(1+\|u\|\big)\Big(u-\frac{v_1(u)}{1+\|u\|}\Big)+\Big(1-\la l_1(u)l_2(u)(1+\|u\|\big)\Big) u\in P_1.
$$
 It follows from the theorem due to Brezis-Martin (see [27]) that $P_1$ is an invariant sets of descending flow of (2.14).  In a similar way we can show that $-P_1$ is also an invariant set of descending flow of (2.14).

Next we  show that $\mbox{int}P_1$  is an  invariant set of (2.14).
 Take  $u\in \mbox{int}P_1$. Note (2.17) also holds  where the integral is in the sense of $E$ topology. Make a variable change $\tau=e^{\mu(s)}-1$ in (2.17). Then we have
$$s=\mu^{-1}\big(\ln(1+\tau)\big),\   ds=\frac{e^{-\mu(s)}}{\mu'(s)}d\tau$$
and
$$\be{l}e^{-\mu(t)}\di\int^t_0e^{\mu(s)}\big(1+\|\si(s, v_0)\|\big)A\si(s, v_0)ds\\=e^{-\mu(t)}\di\int^{e^{\mu(t)}-1}_0 A\si\Big(\mu^{-1}\big(\ln(1+\tau)\big), u\Big)d\tau\ (\mbox{the integral is in the sense of}\ E\ \mbox{topology})\\
=\lim\limits_{n\ri\infty}\frac{1}{n}\big(1-e^{-\mu(t)}\big)\sum\limits_{k=0}^{n-1} A\si\Big(\mu^{-1}\big(\ln(1+\frac{k(e^{\mu(t)}-1)}{n}\big), u\Big).\en $$
It follows from (2.16)  and  Lemma 2.3 that $A(P_1)\subset P_1$.  Since $P_1$ is an invariant set of descending flow of (2.14), we have
$$\bar v:=\lim\limits_{n\ri\infty}\frac{1}{n}\sum\limits_{k=0}^{n-1} A\si\Big(\mu^{-1}\big(\ln(1+\frac{k(e^{\mu(t)}-1)}{n}\big), u\Big)\in P_1.$$
It follows from  (2.17) that
$$\si(t,u)=e^{-\mu(t)} u+(1-e^{-\mu(t)})\bar v.\no(2.19)$$
Since $u\in \mbox{int}P_1$, there exists $r_0>0$ such that $B_1(u, r_0)\subset \mbox{int} P_1$, where $B_1(u, r_0)=\{y\in E:\|y-u\|_1<r_0\}$.
For each $\bar x\in B_1\big(\si(t,u), e^{-\mu(t)} r_0\big)$, let $$\bar y=u+ e^{\mu(t)} (\bar x- \si(t,u)).$$ Then we have
$$\|\bar y-u\|_1=e^{\mu(t)} \|\bar x- \si(t,u)\|_1<r_0.$$
This implies that $\bar y\in B_1(u, r_0)$. On the other hand, by (2.19) we have
$$\bar x=e^{-\mu(t)}\bar y-e^{-\mu(t)} u+\si(t, u)=e^{-\mu(t)} \bar y+(1-e^{-\mu(t)}) \bar v\in P_1.$$
So, $B_1\big(\si(t,u), e^{-\mu(t)} r_0\big)\subset P_1$ and $\si(t,u)\in \mbox{int} P_1$ for all $t\in [0,T_1(v_0))$. This implies that $\mbox{int}\ P_1$ is an invariant set for the descending flow of  (2.14).

Similarly, $\mbox{int}\ (-P_1)$ is an invariant set for the descending flow of  (2.14).

5)\ For each $v_0\in \vr^{-1}\big((-\infty, d_1-\frac{1}{2}]\big)\cap E$ with $\inf\limits_{u\in o(v_0)} \vr(u)\gl-\frac{1}{4}$, we assume the conclusion 5) doesn't hold. It follows from (2.14) and (2.18) that
$$\vr(v_0)-\vr(\si(t, v_0))=-\int^{t}_{0} \frac{\pa}{\pa s} h(s, v_0)ds\gl \frac{\ga}{16} t \ \mbox{for}\  t\in [0,T_1(v_0)).$$
In this case  we have $T_1(v_0)=+\infty$ if $o(v_0)\in D_0\bk K_\de$. So, for $t_0=16\ga^{-1}(d_1+1)$, we have $$\vr(\si(t_0, v_0))\ls \vr(v_0)-\frac{\ga}{16} t_0\ls d_1-\frac{1}{2}-\frac{\ga}{16} t_0<-1,$$
which  contradicts to $\inf\limits_{u\in o(v_0)} \vr(u)\gl-\frac{1}{4}$.

The proof  is complete.
\vskip 0.1in

Similar to the proof of Theorem 1.2 we have the following Lemma 2.5.

{\bf Lemma 2.5.}\ \textit{Let $G\subset E$  be a connected and invariant set of (2.14), and $D$ be an open invariant subset of $G$. Then the following assertions hold:\\
  1)\ $C_G(D)$ is  an open subset of $G$;\\
2)\ $\pa_G C_G(D)$ is  an invariant set of descending flow of (2.14);\\
3)\ $\inf\limits_{u\in \pa_G C_G(D)} \vr(u)\gl \inf\limits_{u\in \pa_G (D)} \vr(u)$.}
\vskip 0.1in

{\bf Lemma 2.6 ([35]).}\  \textit{Assume $U$ is bounded connected open set of $\mathbb{R}^2$ and $(0,0)\in U$,
then there exists a connected component $\Gamma^{\prime}$ of the boundary of $U$, such that each one side ray $l$ emitting from the origin satisfies $l\cap \Gamma^{\prime}\neq \emptyset$.}
\vskip 0.1in

{\bf\textit{ The Proof of Theorem 2.1.}}\
\vskip 0.1in

Take $d_0\in (0, d_1-\frac{1}{2})$. Let $O$ be the connected component of $\vr^{-1}(-\infty, d_0)\cap E$ containing $0$.  It follows from Lemma 2.5 that $C_E(O)$ is an open invariant set of descending flow of (2.14). By Lemma 2.5 (3) and (2.2), we see that $C_E(O)\neq E$, and so $\pa_E C_E(O)\neq\ep$.  By Lemma 2.4 and 2.5 we have
 $$\inf\limits_{u\in \pa_{E} C_{E}(O)} \vr(u)\gl \inf\limits_{u\in \pa_{E} (O)} \vr(u)=d_0.\no(2.20)$$
Since $C_E(O)$ is open in $E$, $C_E(O)\cap E_1\subset B(0, R_0)$ is an open and bounded subset of $E_1$ containing $0$. It follows from Lemma 2.6 that there exists a connected component $\Gamma^{\prime}$ of the boundary of $C_E(O)\cap E_1$, such that each one side ray $l$ emitting from the origin satisfies $l\cap \Gamma^{\prime}\neq \emptyset$. Let $\Gamma$ be the connected component of $\pa_{E} C_{E}(O)$ containing $\Gamma'$. It follows from Lemma 2.5 that $\Gamma$ is an invariant set of descending flow of (2.14).

 It follows from Lemma 2.4 that  $\mbox{int}(P_1)$ is an  invariant set of  descending flow of (2.14). Thus, $\Gamma\cap\ \mbox{int}(P_1)$ is an invariant set of descending flow of (2.14). Take $\wi v_1\in S_{R_0}\cap E_1\cap \mbox{int}(P_1)$ and let $\wi l$ be the ray emitting from the origin and passing through $\wi v_1$. Then, we have $\wi l\cap (\Gamma\cap P_1)\neq\ep$, and so $\Gamma\cap P_1\cap \vr^{-1}\big((-\infty, d_1-1]\big)\neq\ep$. Take $v_0\in \Gamma\cap P_1\cap \vr^{-1}\big((-\infty, d_1-1]\big)$. It follows from (2.20) that $\inf\limits_{u\in o(v_0)} \vr(u)\gl-\frac{1}{4}$, that is $o(v_0)\subset \vr^{-1}\big([d_0, d_1-1]\big)$. Then we have the following two cases:

  1)\ If $v_0\in K_{2\de}$, by (2.3) there must exists a $u_1$ with $u_1\in P_1\cap \vr^{-1}\big([-\frac{1}{4}, d_1]\big)\cap (K\bk\{0\})$.

   2)\ If $v_0\in \big(\Gamma\cap P_1\cap \vr^{-1}([d_0, d_1-1])\big)\bk K_{2\de}$, by (5) in Lemma 2.4 we see that $\si(\tau(v_0), v_0)\in K_{2\de}\cap P_1$ for some $\tau(v_0)>0$, by the same reason as above  we see that there must exists a $u_1$ with $$u_1\in P_1\cap \vr^{-1}\big([-\frac{1}{4}, d_1-\frac{1}{2}]\big)\cap (K\bk\{0\}).$$ Hence, $\vr$ has at least one positive critical point $u_1$.

 Similarly, we can show that $\vr$ has at least one negative critical point $u_2$.

 Now we show that $\vr$ has at least one sign-changing critical point $u_3$. Obviously, $\Gamma\cap \mbox{int}(P_1)$ and $\Gamma\cap \mbox{int}(-P_1)$ are two open invariant sets of descending flow of (2.14) in $\Gamma$. It follows from Lemma 2.5 that $C_\Gamma(\Gamma\cap\mbox{int}(P_1))$ and $C_\Gamma(\Gamma\cap\mbox{int}(-P_1))$ are two open invariant sets of descending flow of (2.14) in $\Gamma$. By the connectedness of $\Gamma$, we see that
 $$D_1:=\Gamma\bk\big(C_\Gamma(\Gamma\cap\mbox{int}(P_1))\cup C_\Gamma(\Gamma\cap\mbox{int}(-P_1))\neq\ep.$$
 Let
 $$D_2:=\Gamma'\bk\big(C_\Gamma(\Gamma\cap\mbox{int}(P_1))\cup C_\Gamma(\Gamma\cap\mbox{int}(-P_1)).$$
Obviously, $D_2\subset D_1$.  Also by the connectedness of $\Gamma'$, we have $D_2\neq\ep$. Take $v_0\in D_2$, by (2.20) we have $\inf\limits_{u\in o(v_0)} \vr(u)\gl-\frac{1}{4}$. In a similar way we can show that $\vr$ has a sign-changing critical point $u_3$. Indeed, if $v_0\in K_{2\de}$, by (2.3) there must exist a $u_3\in (E\bk (P\cup (-P)))\cap \vr^{-1}([-\frac{1}{4}, d_1-\frac{1}{2}])\cap K$. If $v_0\in \big(D_2\cap \vr^{-1}([d_0, d_1-\frac{1}{2}])\big)\bk K_{2\de}$, by (5) in Lemma 2.4 we see that $\si(\tau(v_0), v_0)\in K_{2\de}\cap D_1$ for some $\tau(v_0)>0$, by the same reason as above  we see that there must exists a $u_3$ with $$u_3\in (E\bk (P\cup (-P)))\cap \vr^{-1}([-\frac{1}{4}, d_1-\frac{1}{2}])\cap (K\bk\{0\}).$$ Hence, $\vr$ has at least one sign-changing critical point $u_3$.  The proof is complete.
 \vskip 0.1in

 {\bf Remark 2.3.}\ It is important in the proof that  $P_1$ and $-P_1$ have nonempty interiors. Instead of these two sets,   the authors of  [28] used the neighborhoods of cones and negative cones. Therefore, we can establish a result similar to Theorem 2.1  based on that condition of [28]. However, for the sake of brevity, we will not discuss the relevant results in this article.
 \vskip 0.1in

Assume that $(I-\na \vr)(\pm P)\subset \pm P$. As in the proof of  Lemma 2.3, for each $x_0\in P$,  there exists $(u_2(x_0), x_0^*)\in \big(\{x_0\}-P)\cap \bar B(0,1)\big)\ti \pa\vr(x_0)$ satisfying  $m_P(x_0)=\lge x_0^*, u_2(x_0))$.  Take  $y_0\in \{x_0\}-\textsf F(x_0^*)$. It follows from the condition $(I-\na \vr)( P)\subset P$ that $y_0\in P$.  So, if $\|x_0-y_0\|<1$,
$$m_P(x_0)=\lge x_0^*, u_2(x_0))\gl \lge x_0^*, x_0-y_0\rge=\lge x_0^*, \textsf F(x_0^*)\rge=\|x_0^*\|_*^2\gl m^2(x_0); $$
if $\|x_0-y_0\|\gl 1$, let $z_0=x_0+\frac{y_0-x_0}{2\|x_0-y_0\|}$,  we have $z_0\in P$ and $\|x_0-z_0\|=\frac{1}{2}$, and so
$$\be{ll}m_P(x_0)&=\lge x_0^*, u_2(x_0))\gl \lge x_0^*, x_0-z_0\rge\\
&=\frac{1}{2\|x_0-y_0\|}\lge x_0^*,x_0-y_0\rge\\
&=\frac{1}{2\|x_0\|_*}\|x_0^*\|_*^2=\frac{1}{2}\|x_0^*\|_*\gl\frac{1}{2} m(x_0).\en$$
Thus, we have $$m_P(x)\gl \min\{\frac{1}{2}, m(x)\}m(x)\ \mbox{for all}\  x\in P.\no(2.21)$$
Similarly, we have
$$m_{-P}(x)\gl \min\{\frac{1}{2}, m(x)\}m(x)\ \mbox{for all}\  x\in -P.\no(2.22)$$
\vskip 0.1in

{\bf Definition 2.6.}\ \textit{We say that $\vr$ satisfies the $(PS)$ condition, for every sequence $\{x_n\}\subset X$ such that $\{\vr(x_n)\}$ is bounded and $ m(x_n)\ri 0$ as $n\ri \infty$, has a convergent subsequence.}
\vskip 0.1in

By using (2.21) and (2.22), in a similar way to show  Theorem 2.1 we can prove  the following Theorem 2.2.

{\bf Theorem  2.2.}\  \textit{Suppose that  (H$_2)$, (H$_3)$, (H$_5)$ and (H$_6)$ hold, $(I-\na \vr)(\pm P)\subset \pm P$, $\vr$ satisfies the condition $(PS)$. Then $\vr$ has at least one positive, one negative and one sign-changing critical point.}
\section{Applications to  Differential Inclusion Problems}

\ \ \
Consider the following  differential  inclusions problems
$$\left\{\be{ll} -\mbox{div}\big(\|Du(x)\|^{p-2} Du(x)\big)-\la |u(x)|^{p-2} u(x) \in  \pa j(x, u) &\ \mbox{in}\ \Om,\\
u=0\ &\ \mbox{on}\ \pa\Om,\en\right.\no(3.1_\la)
$$
where $\Om $ ia a bounded open domain in $\mb R^N$ with  a $C^{1,\al}$-boundary $\pa \Om$ $(0<\al<1)$, $1<p<+\infty$,  the reaction term $\pa j(x,s)$ is the generalized gradient of a non-smooth potential $s\mapsto j(x,s)$, which is subject to the following conditions.

{\bf (H$_j)$}\ $j:\Om\ti \mb R \mapsto \mb R$ is a Carath\'{e}odory  function and there exist constants $a_1>0$, $p<q<p^*$ such that

(i)\ $j(x,\cdot)$ is locally Lipschitz for almost every $x\in \Om$ and $j(x, 0)=0$ a.e. on $\Om$;

(ii)\ $|\xi|\ls a_1(1+|s|^{q-1})$ a.e. in $\Om$ and  for all $s\in \mb R$, $\xi\in\pa j(x,s)$;

(iii)\ there exist constants $\mu>p$ and $M>0$ such that $$\inf\limits_{x\in \Om} j(x, M)>0\ \mbox{and}\ \mu j(x, z)\ls-j^o(x, z;-z)\ \mbox{a.e. on}\ x\in \Om\ \mbox{all}\ z\gl M;$$

(iv)\ $$\lim\limits_{z\ri 0}\sup\frac{pj(x,z)}{z^p}= 0$$ uniformly with respect to $x\in \Om$;

(v)\ $zw(x)\gl 0$ for each $w(x)\in \pa j(x, z)$ a.e. $\Om$ and $z\in\mb R$.
\vskip 0.1in

{\bf Remark 3.1.}\ The main purpose of this part is to show that our theoretical results can be applied to the study of differential inclusion problems. It should be pointed out that some of the conditions we listed above are similar to those in [18], and the proofs of some of the following lemmas are also similar to those in [18,23].
\vskip 0.1in

Let  $$\|u\|=\Big(\int_\Om\|D u\|^p dx\Big)^{1\over p}, |u|_r=\Big(\int_\Om| u|^r dx\Big)^{1\over r}$$
be the standard norms of  $W_0^{1,p}(\Om)$, respectively   $L^r(\Om)$ for $1<r<p^*$.
Let $X=W_0^{1,p}(\Om)$ and $E=C_0^1(\Om)$. For $\la>0$, we introduce the energy functional $\vr_\la: X\mapsto \mb R$  by
 $$\vr_\la (u)=\frac{1}{p}\|u\|^p_p-\frac{\la}{p}|u|_p^p-\int_\Om j(x, u(x))dx,$$
Let $P=\{u\in X: u(x)\gl 0\  \mbox{a.e.}\ \Om\}$ and $P_1=P\cap E$. Given $\la>0$, we say that $u\in X$ is a (weak) solution of (3.1$_\la$) if $\Delta_p u\in L^{q'}(\Om)$, where $\frac{1}{q}+\frac{1}{ q'}=1$, and
$$-\Delta_p u(x)\in\la |u(x)|^{p-2}u(x)+ \pa j(x, u(x))\ \ \mbox{for almost every}\ x\in \Om.$$
Let $K_\la=\{x\in X: 0\in \pa\vr_\la(x)\}$ for $\la>0$.
 \vskip 0.1in

Recall some facts about the spectrum of the  $p$-Laplacian with  Dirichlet boundary  condition. Consider the nonlinear eigenvalue problem
$$\left\{\be{ll} -\mbox{div}\big(\|Du(x)\|^{p-2} Du(x)\big)=\la |u(x)|^{p-2} u(x) &\ \mbox{in}\ \Om,\\
u=0\ &\ \mbox{on}\ \pa\Om,\en\right.\no(3.2_\la)
$$
Let $\la_1$ be the principal eigenvalue of $(-\Delta_p, W_0^{1,p}(\Om))$. Then $\la_1$ is positive, isolated and simple. There is the following variational characterization of $\la_1$ using Rayleigh quotient:
$$\la_1=\inf\Big\{\frac{\|Du\|_p^p}{|u|_p^p}: u\in W_0^{1,p}(\Om), u\neq 0\Big\}.$$
This minimum is actually realized at normalized eigenfunction $u_1$. The Ljusternik-Schnirelmann theory gives, in addition to $\la_1$, a whole strictly increasing sequence of positive numbers $\la_1<\la_2\ls\la_3\ls\cdots\ls \la_k\ls\cdots$ for which there exist nontrivial solutions for problem (3.2$_\la$). In what follows we let $u_2\in W_0^{1,p}(\Om)$ be a nontrivial solutions for problem (3.2$_\la$) corresponding to $\la_2$, and $E_1=\mbox{span}\ \{u_1, u_2\}$.
\vskip 0.1in

{\bf Theorem 3.1.}\  \textit{Assume (H$_j)$ holds. Then for $\la\in (0,\la_1)$, (3.1$_\la$) has at least one positive solution, one negative solution and one  sign-changing solution.}
\vskip 0.1in

{\bf Lemma 3.1.}\ \textit{If $u\in K_\la$, then $u\in C_0^1(\bar\Omega)$ and $u$ solves (3.1$_\la)$. Moreover, if $u\in \pm P\cap K_\la$ and $u\neq 0$, then $u\in \mbox{int} (\pm P_1)\cap K_\la$. }

{\bf Proof.}\ The proof is similar to Proposition 3.1 and 3.2 in [23]. Let $A: W_0^{1,p}(\Om)\mapsto W^{-1, p'}(\Om)$ be defined by
$$\lge A(u), v\rge=\int_\Om \|D u(z)\|^{p-2}(Du(z), Dv(z))_{\mb R^N} dz \ \ \mbox{for}\ u,v\in W_0^{1,p}(\Om).$$
It is known that $A$ is monotone and demi-continuous, hence maximal monotone, and so generalized pseudomonotone.

 Obviously, $u\mapsto \frac{1}{p} \|u\|^p$ is a $C^1$-functional whose derivative is the   operator $A$. Aubin-Clarke's Theorem ensures  that the functional $$
u\mapsto \int_\Omega j(x, u)dx$$
is Lipschitz continuous on any bounded subset of $L^q(\Om)$ and its gradient is included in the set
$$N(u)=\{ w\in L^{q'}(\Om): w(x)\in\pa j(x, u(x))\ \mbox{for almost every}\ x\in\Om\}.$$
Since $X$ continuously embedded in $L^q(\Om)$, the function $\vr_\la$ turns out to be locally Lipschitz on $X$. So, we have
$$ \pa \vr_\la (u)\subset A(u)-\la|u|^{p-2} u- N(u).\no(3.3)$$
Now, if $u\in X$  complies with $0\in \pa \vr_\la(u)$ then $$A(u)=\la |u|^{p-2}u+w\ \mbox{in}\ X^*$$
for some $w\in N(u)$. Hence, $\Delta_p u\in L^{q'}(\Omega)$ and $u$ solves (3.1$_\la)$. By the condition (H$_j)$ (ii) and (3.3) we get the estimate
$$-u\Delta_p u\ls a_1 (|u|+|u|^q)\ \mbox{a.e. in}\  \Omega.$$
Hence, by [12, Theorem 1.5.5], we have $u\in L^\infty(\Omega)$.  From (H$_j)$ (ii) it follows $\Delta_p u\in L^\infty(\Omega)$. So, by [12, Theorem 1.5.6], we have $u\in C_0^1(\bar \Omega)$.

Let $u\in P\cap K$ and $u\neq 0$. By (H$_j)$ (v), for each  $c_0>0$ we have
$$\Delta_p u=-\la u^{p-1}-w\ls c_0 u^{p-1}$$
for some $w\in \pa j(x, u)$. The V\'{a}zquez maximum principle yields $u\in\mbox{int}(P_1)$.

Similarly, if $u\in -P\cap K_\la$ and $u\neq 0$,  then $u\in \mbox{int} (-P_1)\cap K_\la$. The proof is complete.
\vskip 0.1in

{\bf Lemma 3.2.}\ \textit{The functional $\vr_\la: X\mapsto\mb R$ satisfies the conditions  $(CPS)$, $(CPS)_P$ and $(CPS)_{-P}$ for $\la>0$.}

{\bf Proof.}\ The proof is similar to claim 1 of Theorem 3.1 in [18]. For reader's convenience we give the details of the process. We divide the proof into three steps.

{\bf Step 1.}\  In this step we prove the assertion: \textit{if $\{x_n\}\subset W_0^{1,p}(\Om)$ is bounded, and either $(1+\|x_n\|)m(x_n)\ri 0$, or $(1+\|x_n\|)m_{\pm P}(x_n)\ri 0$ as $n\ri +\infty$, then $\{x_n\}$ has a convergent subsequence.}

 We only consider the case of  $(1+\|x_n\|)m(x_n)\ri 0$ as $n\ri +\infty$. In a similar way we can prove the case of $(1+\|x_n\|)m_{\pm P}(x_n)\ri 0$.

 Since $\{x_n\}\subset W_0^{1,p}(\Om)$ is bounded, by passing to a subsequence if necessary, we may assume
$$x_n\rightharpoonup x\ \mbox{in}\ W_0^{1,p}(\Om), x_n \ri x \ \mbox{in}\ L^r(\Om)\ \mbox{for}\ 1<r<p^*, x_n(z)\ri x(z)\ \mbox{a.e. on }\ \Om$$
and $|x_n(z)|\ls k(z)\ \mbox{a.e. on }\ \Om$, for all $n\gl 1$, with $k\in L^{q'}(\Om)$.  Take $x_n^*\in \pa \vr_\la(x_n)$ such that $m(x_n)=\|x_n^*\|_{*}$ for $n\gl 1$. Then we have $$x_n^*=A (x_n)-\la |x_n|^{p-2} x_n- u_n\no(3.4)$$
with $u_n\in L^{q'}(\Om)$, satisfying  $u_n(x)\in \pa \vr(x, x_n(x))$ a.e. on $\Om$.

Now, we can deduce from $(1+\|x_n\|)m(x_n)\ri 0$ that   $|\lge x_n^*, x_n-x\rge |\ls \frac{1}{n}\|x_n-x\|$. This reads
$$\Big|\lge A(x_n), x_n-x\rge-\la\int_\Om |x_n|^{p-2} x_n (x_n-x) dz-\int_\Om u_n(x_n-x)dz\Big|\ls\frac{1}{n}\|x_n-x\|.$$
Consequently, we have $$\la\int_\Om |x_n|^{p-2} x_n (x_n-x) dz\ri 0\ \mbox{and}\ \int_\Om u_n(x_n-x)dz\ri 0\ \mbox{as} \ n\ri\infty, $$
and so, $$\lim\limits_{n\ri\infty}\lge A(x_n), x_n-x\rge=0.$$
Since $A$ is  generalized pseudomonotone, we have  $\lge A(x_n),x_n\rge\ri \lge A(x), x\rge$, or equivalently, $\|Dx_n\|_p\ri\|Dx\|_p$. Recalling that $Dx_n\rightharpoonup Dx$ in $L^p(\Om, \mb R^N)$  and  $L^p(\Om, \mb R^N)$ being uniformly convex we have  $Dx_n\ri Dx$ in $L^p(\Om, \mb R^N)$ which  means $x_n\ri x$ in $W_0^{1,p}(\Om)$.

{\bf Step 2.}\  $\vr_\la$ satisfies the condition $(CPS)$.

Let $\{ x_n\}\subset X$  be such that $|\vr_\la(x_n)|\ls M_1$ for some $M_1>0$ and $(1+\|x_n\|)m(x_n)\ri 0$ as $n\ri \infty$, and so $m(x_n)\ri 0$ as $n\ri\infty$.
Take $x_n^*\in \pa \vr(x_n)$ such that $m(x_n)=\|x_n^*\|_{*}$ for $n\gl 1$. Then (3.4) holds.  Since $m(x_n)\ri 0$ as $n\ri \infty$, we can say that $|\lge x_n^*, x_n\rge |\ls \frac{1}{n}\|x_n\|$. So,
$$-\|Dx_n\|_p^p+\la|x_n|_p^p-\int_\Om j^o(z, x_n(z); -x_n(z))dz\ls \frac{1}{n}\|x_n\|.\no(3.5)$$
Similarly, since $|\vr_\la (x_n)|\ls M_1$ for all $n\gl 1$, we have
$$\frac{1}{p}\|Dx_n\|_p^p-\frac{\la}{p}|x_n|_p^p-\int_\Om j(z, x_n(z))dz\ls M_1.\no(3.6)$$
By  (3.5) and (3.6) we obtain
$$\be{l} \Big(\frac{\mu}{p}-1\Big)\|D x_n\|_p^p-\la\big(\frac{\mu}{p}-1\big)|x_n|_p^p\\
-\di\int_\Om \big( \mu j(z, x_n(z))+ j^o(z, x_n(z); - x_n(z))\big)dz\ls \mu M_1+\frac{1}{n}\|x_n\|.
\en$$
By (H$_j)$ (ii)(iii) we have  for some $\beta_1>0$,
$$\di\int_\Om \big( \mu j(z, x_n(z))+ j^o(z, x_n(z); - x_n(z))\big)dz\gl -\beta_1.$$
Then, we have for $n\gl 1$,
$$\Big(\frac{\mu}{p}-1\Big)\Big(1-\frac{\la}{\la_1}\Big)\|Dx_n\|_p^p\ls \mu M_1+\frac{1}{n}\|x_n\|+\beta_1.$$
Because $\la<\la_1$ and $\mu> p$,  we infer that $\{x_n\}\subset W_0^{1,p}(\Om)$ is bounded. So, by the assertion in  Step 1 we see that $\{x_n\}$ has a convergent subsequence.

{\bf Step 3.}\  $\vr_\la$ satisfies the condition $(CPS)_P$ and $(CPS)_{-P}$.

 Let $\{ x_n\}\subset P$  be such that $|\vr_\la(x_n)|\ls M_2$ for some $M_2>0$,  and $(1+\|x_n\|)m_{P}(x_n)\ri 0$ as $n\ri \infty$.
Let $h_p: X^*_w\ti X\mapsto \mb R$ be defined by
$$h_{P}(x^*, x)=\sup\big\{\lge x^*, x-y\rge: y\in P, \|y-x\|<1\big\}.$$
Then, $h_p: X^*_w\ti X\mapsto \mb R$  is lower semi-continuous; See [Lemma 1, 14] or [Lemma 2.2.1, 12]. Now we can find  $x_n^*\in \pa \vr_\la(x_n)$ such that $m_{P}(x_n)=h_P(x_n^*, C_1(x_n))$, where $C_1(x_n)=(x_n-P)\cap B(0,1)$. So, we have
$$\lge x_n^*, x_n-y\rge\ls h_P(x_n^*, C_1(x_n))\ \mbox{for all}\ y\in P, \|x_n-y\|<1.\no(3.7)$$
Let $y_n=x_n+\frac{1}{2\|x_n\|} x_n$ for $n\gl 1$. Then, $y_n\in P$, $\|x_n-y_n\|<1$. And so,  by (3.7) we have
$$\be{ll} (1+\|x_n\|)\lge x_n^*, x_n- y_n\rge&=-\frac{(1+\|x_n\|)}{2\|x_n\|}\lge x_n^*,x_n\rge\\
&=-\frac{(1+\|x_n\|)}{2\|x_n\|}\Big(\lge A(x_n),x_n\rge-\la|x_n|_p^p+\di\int_\Om j^o(z, x_n(z); -x_n(z))dz\Big)\\
&\ls (1+\|x_n\|)h_P(x_n^*, C_1(x_n))=:\va_n\ \mbox{with}\ \va_n\downarrow 0.\en
$$
So we obtain $$-\Big(\lge A(x_n),x_n\rge-\la|x_n|_p^p+\di\int_\Om j^o(z, x_n(z); -x_n(z))dz\Big)\ls \frac{2\|x_n\|}{(1+\|x_n\|)}\va_n.$$
which reads
$$-\|Dx_n\|_p^p+\la |x_n|_p^p+\di\int_\Om j^o(z, x_n(z); -x_n(z))dz \ls \frac{2\|x_n\|}{(1+\|x_n\|)}\va_n=\va'_n, \va_n'\downarrow 0.$$
Then, by using the arguments from (3.5) to the end of the step 2 we can show that $\{\|x_n\|\}$ is bounded in $W_0^{1,p}(\Om)$. By the assertion in Step 1 we see that $\{x_n\}$ has a convergent subsequence. Hence, $\vr_\la$ satisfies the condition $(CPS)_P$.

Similarly, we can show that $\vr_\la$ satisfies the condition  $(CPS)_{-P}$.
The proof is complete.
\vskip 0.1in

{\bf Lemma 3.3.}\ \textit{There exists $R_0>0$ such that for $\la\in (0,\la_1)$,
 $$\sup\limits_{u\in S_{R_0}\cap E_1} \vr_\la(u)< 0.\no(3.8)$$}

{\bf Proof.}\ The proof is similar to claim 3 of Theorem 3.1 in [18].  For almost all $x\in \Om$ and all $z\in \mb R$, the function $s\mapsto \frac{1}{s^\mu} j(x, sz)$ is locally Lipschitz on $(0, +\infty)$. Using the mean value theorem for locally Lipschitz functions, for $s>1$ we can find $\theta\in (1,s)$ such that
$$\be{ll}\frac{1}{s^\mu}j(x, sz)-j(x, z)&\in \Big(-\frac{\mu}{\theta^{\mu+1}}j(x, \theta z)+\frac{1}{\theta^\mu}\pa_z(x, \theta z) z\Big)(s-1)\\
&=\frac{s-1}{\theta^{\mu+1}}\big(-\mu j(x, \theta z)+\pa_z j(x, \theta z)\theta z\big).\en
$$
 By (H$_j)$ (iv), for almost all $x\in \Om$ and all $z\gl M$,  we have
 $$\frac{1}{s^\mu}j(x, sz)-j(x, z)\gl \frac{s-1}{\theta^{\mu+1}}\big(-\mu j(x, \theta z)-j^o(x, \theta z; -\theta z)\big)\gl 0.$$
 Then   for almost all $x\in \Om$ and all $z\gl M$,  we have
 $$j(x,z)=j\big(x, \frac{z}{M}M\big)\gl \big(\frac{z}{M}\big)^\mu j(x, M)\gl \big(\frac{z}{M}\big)^\mu\inf\limits_{x\in\Om} j(x, M).$$
 Let $p_1\in (p, \mu)$.  So, it is seen that   for a given $\eta>0$ we can find a constant $c_\eta>0$ such that
 $$j(x,z)\gl \frac{\eta}{p} z^{p_1}-c_\eta\ \mbox{for}\ \mbox{a.e.}\ x\in\Om.\no(3.9)$$
 Let $u_0\in S_1:=\{u\in W_0^{1, p}(\Om): \|Du\|_p=1\}$. It follows from (3.9) that
 $$\be{ll}\vr_\la(t u_0)&=\frac{1}{p}\|D(tu_0)\|_p^p-\frac{\la t^p}{p}|u_0|_p^p-\di\int_\Om j\big(x, tu_0(x)\big) dx\\
 &\ls \frac{t^p}{p}-\frac{\la t^p}{p}|u_0|_p^p-\frac{\eta t^{p_1}}{p} |u_0|_p^p+c_\eta|\Om|.\en $$
 This implies that
$$\lim\limits_{t\ri+\infty}\vr_\la(tu_0)=-\infty. $$
 Since $E_1\cap S_1$ is compact, there exists $R_0>0$ such that (3.8) holds. The proof is complete.
\vskip 0.1in

{\bf Lemma 3.4. }\ \textit{Let $\la\in (0,\la_1)$. There exists  $R_1\in (0, R_0)$ such that  $\beta_d:=\inf\limits_{u\in S_d} \vr_\la(u)>0$ for any $d\in (0, R_1)$.}

{\bf Proof.}\  By the condition (H$_j)$ (ii) (iv), for $\va>0$ with $\la+\va<\la_1$, we can find  a $c_\va>0$ such that
 $$0\ls j(x, z)\ls \frac{\va}{p}|z|^p+c_\va |x|^q\  \mbox{for a.e.}\  x\in \Om, z\in \mb R.$$
So, we have
$$\vr_\la (u)\gl \frac{1}{p} \Big(1-\frac{\la+\va}{\la_1}\Big)\|Du\|_p^p- c_1\|Du\|_p^q,\ \  \forall u\in W_0^{1,p}(\Om)$$
for some $c_1>0$. Thus, there exists $R_1\in (0, R_0)$ small enough such that $\beta_d:=\inf\limits_{u\in S_d} \vr_\la(u)>0$ for any $d\in (0, R_1)$, and $0$ is a strictly local minimum point. The proof is complete.
\vskip 0.1in

{\bf Lemma 3.5. }\ \textit{$\vr_\la$ is outwardly directed on $\pm P$ for $\la\in (0,\la_1)$.}

{\bf Proof.}\ The  following elementary inequality is well known:
$$\big(|y|^{p-2} y-|h|^{p-2}h, y-h\big)_{\mb R^N}\gl\left\{\be{ll}c_1(p)(|y|+|h|)^{p-2}|y-h|^2\ &\mbox{if}\ 1<p<2,\\
c_2(p)|y-h|^p  &\mbox{if}\ p\gl 2\en\right.
$$
for all $y,h\in\mb R^N$, where $c_1(p), c_2(p)>0$ are constants.

For each $u\in P$, $u^*=Au-\la |u|^{p-2} u-w\neq 0$, let $v=(-\Delta_p)^{-1}\big(\la |u|^{p-2} u+w)$, where $w\in L^{q'}(\Om)$ and $w\in \pa j(x, u)$. Then, we have
$$\be{ll}\lge u^*, u-v\rge&=\lge Au-\la |u|^{p-2} u-w, u-v\rge\\
&=\lge Au+\Delta_p v, u-v\rge =\lge Au, u-v\rge+\lge \Delta_p v, u-v\rge\\
&=\di\int_\Om|\na u|^{p-2}\na u\cdot\na (u-v)+\di\int_\Om(u-v)\Delta_p v\\
&=\di\int_\Om|\na u|^{p-2}\na u\cdot\na (u-v)-\di\int_\Om |\na v|^{p-2}\na v\cdot\na (u-v)\\
&\gl \left\{\be{ll}c_1(p)\di\int_\Om (|\na u|+|\na v|)^{p-2}|\na (u-v)|^2\ &\mbox{if}\ 1<p<2,\\
c_2(p)\di\int_\Om |\na(u-v)|^p  &\mbox{if}\ p\gl 2\en\right.\\
&> 0.\en
$$
 This implies that $\vr_\la$ is outwardly directed on $P$. Similarly, $\vr_\la$ is outwardly directed on $-P$. The proof is complete.
\vskip 0.1in

{\bf Lemma 3.6.}\ \textit{ For any $a,b\in\mb R$ with $a<b$, $K_\la\cap \vr_\la^{-1}([a,b])$ is compact in $E$.}

{\bf Proof.}\ For each $u\in K_\la\cap \vr_\la^{-1}([a,b])$, it follows from the proof of Lemma 3.1 that
$u=(-\De_p)^{-1}(\la |u|^{p-2}u+w)$ for some $w\in N(u)$.  By a similar way as the proof of Lemma 3.2 we can prove that  $K_\la\cap \vr_\la^{-1}([a,b])$ is bounded in $X$.  Let
$$B(\la):=\big\{\la |u|^{p-2}u+w: u\in K_\la\cap \vr_\la^{-1}([a,b]), w\in N(u)\big\}.$$

If $p>N$ then $X\hookrightarrow L^\infty(\Om)$. So, $K_\la\cap \vr_\la^{-1}([a,b])$ is bounded in $L^\infty(\Om)$ if $p>N$. It follows the condition (H$_j)$ (ii)  that  the set $B(\la)$
is bounded in $L^\infty(\Om)$. According to [33], there exists $0<\al<1$  and $c_1>0$ such that
$$\|(-\De_p)^{-1}u\|_{C^{1,\al}}\ls \|u\|_{\infty}^{\frac{1}{p-1}}\ \  \mbox{for all}\ u\in L^\infty(\Om).\no(3.10)$$
Hence, $K_\la\cap \vr_\la^{-1}([a,b])$ is bounded in $C^{1,\al}(\bar\Om)$, and is compact in $E$.

If $1<p\ls N$, take $p^*>r>\frac{q N}{p}$. It follows the condition (H$_j)$ (ii)  that  the set $B(\la)$
is bounded in $L^r(\Om)$. According to  [34] we have
$$\|(-\De_p)^{-1}u\|_\infty\ls c_2\|u\|_r^{\frac{1}{p-1}}\ \ \mbox{for all}\ u\in L^r(\Om).\no(3.11)$$
It follows from  (3.11)  that $K_\la\cap \vr_\la^{-1}([a,b])$ is bounded in $L^\infty(\Om)$. Then, by  (3.10), we see that $K_\la\cap \vr_\la^{-1}([a,b])$ is compact in $E$. The proof is complete.
\vskip 0.1in

{\bf\textit{ The Proof of Theorem 3.1.}}\ The conclusion can be easily proved by Lemma 3.1$\sim$3.6, and Theorem 2.1. The proof is complete.
\vskip 0.1in

{\bf Remark 3.2.}\ There have been some papers studied the existence for sign-changing solutions of differential inclusion problems; see [23,24] and the references therein. For example, by combing variational methods with truncation techniques the paper [23] obtained the existence of positive, negative and nodal solutions to differential inclusion problems with a parameter. Here, our  method is different to  that in [23,24].
\vskip 0.1in

{\bf Remark 3.3.}\ Zhang, Chen and Li [7],   and Bartsch T. and Liu [8]  study  $p$-Laplacian  equations boundary value problems with  smooth nonlinearities. They proved some results of at least three solutions: one positive, one negative, and one sign-changing. To show their main results, they  construct of special continuously differentiable  pseudo-gradient vector fields which ensuring that both the cone and negative cone are  invariant set of descending flow generated by that pseudo-gradient vector fields. Theorem 3.1 can be thought as an extension of some of the main results in [7,8]. Here, we allow the nonlinearity to have discontinuity. Moreover, our method to construct the pseudo-gradient vector fields is different to that in [7,8].

\end{document}